\documentclass[mmnp]{my-edpsmath}

\usepackage{graphicx}
\usepackage{subfig}
\usepackage{url}
\usepackage{amsmath,amssymb,amsfonts,amsthm}
\usepackage{hyperref}
\usepackage{enumerate}

\begin{document}
	
\title{A SIQRB delayed model for cholera and optimal control treatment}

\runningtitle{A SIQRB delayed model for cholera and optimal control treatment}

\thanks{This research was supported by the
Portuguese Foundation for Science and Technology (FCT)
within projects UIDB/04106/2020 and UIDP/04106/2020 (CIDMA). 
Lemos-Pai\~{a}o was also supported by the Ph.D.
fellowship PD/BD/114184/2016; Silva 
by FCT via the FCT Researcher Program CEEC 
Individual 2018 with reference CEECIND/00564/2018.}
\thanks{Corresponding author: cjoaosilva@ua.pt}


\author{Ana P. Lemos-Pai\~{a}o}
\address{Center for Research and Development in Mathematics and Applications (CIDMA), 
Department of Mathematics, University of Aveiro, 3810-193 Aveiro, Portugal;
\email{anapaiao@ua.pt\ \&\ cjoaosilva@ua.pt\ \&\ delfim@ua.pt}}

\author{Helmut Maurer}
\address{Institute of Computational and Applied Mathematics,
University of M\"{u}nster, D-48149 M\"{u}nster, Germany;\\ 
\email{helmut.maurer@uni-muenster.de}}

\author{Cristiana J. Silva}
\sameaddress{1}

\author{Delfim~F.~M.~Torres}
\sameaddress{1}

\runningauthors{A. P. Lemos-Pai\~{a}o, H. Maurer, C. J. Silva, D. F. M. Torres}

\date{Submitted to MMNP: Sept 23, 2021; Revised: March 15, 2022; Accepted: June 25, 2022} 


\begin{abstract} 
We improve a recent mathematical model for cholera  
by adding a time delay that represents the time between the instant at which an individual 
becomes infected and the instant at which he begins to have symptoms of cholera disease.
We prove that the delayed cholera model is biologically meaningful and analyze the local 
asymptotic stability of the equilibrium points for positive time delays.
An optimal control problem is proposed and analyzed, where the goal is to obtain optimal 
treatment strategies, through quarantine, that minimize the number of infective individuals 
and the bacterial concentration, as well as treatment costs. Necessary optimality 
conditions are applied to the delayed optimal control problem, with a $L^1$ type cost functional.
We show that the delayed cholera model fits better the cholera outbreak 
that occurred in the Department of Artibonite -- Haiti, from 1 November 2010 
to 1 May 2011, than the non-delayed model. Considering the data 
of the cholera outbreak in Haiti, we solve numerically the delayed 
optimal control problem and propose solutions 
for the outbreak control and eradication.
\end{abstract}

\subjclass{34C60, 49K15, 92D30.}

\keywords{SIQRB cholera delayed model;
time-delay;
disease-free and endemic equilibria;
stability;
optimal control;
case study in Haiti.}

\maketitle


\section{Introduction}

Cholera is an acute diarrheal infectious disease caused by infection of
the intestine with the bacterium \emph{Vibrio cholerae}, which lives in an aquatic organism.
There are 200 serogroups of the bacterium \emph{Vibrio cholerae}, but only two of them 
(O1 and O139) are responsible for the cholera disease \cite{Cui,Kirschner}. 
They pass through and survive the gastric acid barrier of the stomach and then 
penetrate the mucus lining that coats the intestinal epithelial \cite{Cui,Reidl}. 
They colonize the intestine and then produce enterotoxin, 
which stimulates water and electrolyte secretion by the endothelial 
cells of the small intestine \cite{Cui}. The ingestion of contaminated 
food or water can cause cholera outbreaks,
as proved by John Snow in 1854 \cite{Shuai}.
Nevertheless, there are other ways of spreading.
Susceptible individuals can also become infected
if they contact with infectious individuals. If these individuals
are at an increased risk of infection, they can transmit
the disease to other people that live with them
and are involved in food preparation or use water storage containers, see e.g. \cite{Shuai}.
An individual can be infective with or without symptoms which can appear from 
a few hours to 5 days, after infection. However, symptoms typically appear in 2--3 days: 
see \cite{TimeDelayCholera}. Some symptoms are vomiting, leg cramps and copious, 
painless, and watery diarrhea. It is very important that infective individuals 
can get treatment as soon as possible, because without it they become dehydrated 
and they can suffer from acidosis and have a circulatory collapse. Even worse, 
this situation can lead to death within 12 to 24 hours \cite{Shuai,Mwasa}.
Some studies and experiments suggest that a recovered individual can be immune
to the disease during a period of 3 to 10 years. On the other hand, recent 
researches suggest that immunity can be lost after a period of weeks to months \cite{Shuai,Neilan}.
Diseases involving diarrhea are the major cause of child mortality in developing countries, 
because the access to clean drinking water and sanitation is difficult \cite{Beryl}. 
Moreover, Sun et al. wrote in \cite{Sun} that this disease has generated a great threat 
to human society and caused enormous morbidity and mortality with weak surveillance system. 
Thus, it is very important to study mathematical models of the cholera spread 
in order to know how to curtail it.

Several mathematical models for the spread of cholera have been studied since, at least, 1979: 
see, e.g., \cite{Shuai,Mwasa,Neilan,Capasso:1979,Capone,Codeco:2001,Hartley:2006,%
Hove-Musekwa,Joh:2009,Lemos-Paiao,Mukandavire:2008,Nishiura:2018,Nishiura:2017,Pascual,Wang,Liu} 
and references cited therein. In \cite{Neilan}, a SIR (Susceptible--Infectious--Recovered)
type model is proposed, which considers two classes for the bacterial concentrations
(hyper-infectious and less-infectious) and two classes for the infective individuals 
(asymptomatic and symptomatic). The authors compare a cost-effective balance of multiple 
intervention methods of two endemic populations, using optimal control theory, parameter 
sensitivity analysis, and numerical simulations. Wang and Modnak \cite{Wang} also consider 
a SIR type model with a class for the \emph{Vibrio cholerae} concentration in the environment. 
The model incorporates three control measures: vaccination, therapeutic treatment, 
and water sanitation. The stability analysis of equilibrium points is done when the controls 
are given by constant values. They also study a more general cholera model with time-dependent 
controls, proving existence of solution to an optimal control problem and deriving 
necessary optimality conditions based on Pontryagin's Maximum Principle.
The authors of \cite{Mwasa} incorporate in a SIR type model
public health educational campaigns, vaccination, quarantine, and treatment,
as control strategies. The model also considers a class for the bacterial concentration.
The education-induced, vaccination-induced, and treatment-induced reproductive numbers, 
as well as the combined reproductive number, are compared with the basic reproduction number 
to assess the possible community benefits of these strategies. The stability analysis 
of the equilibria is performed using a Lyapunov functional approach.
In \cite{Shuai}, a SIR type model with distributed delays is proposed. It incorporates
hyperinfectivity, where infectivity varies with the time since the pathogen was shed,
and temporary immunity. The basic reproduction number is computed and it plays 
an important role to know if the disease dies out or not. Numerical simulations 
are carried out in order to illustrate important details of the unique endemic equilibrium's stability.
In \cite{Wang2}, Wang and Liao present a SIR type model with a class for the 
\emph{Vibrio cholerae} concentration in the contaminated water. It is an unified 
deterministic model for cholera, because it considers a general incidence rate 
and a general formulation of the pathogen concentration. The basic reproduction 
number is computed and conditions are derived for the existence of the disease-free 
and endemic equilibrium points. The local asymptotic and global stability analysis 
of the equilibrium points are studied. The authors show that different models 
can be studied in a single framework, using three representative cholera 
models presented in \cite{Codeco:2001, Hartley:2006,Mukandavire:2011}.
A mathematical model that considers public health educational campaigns, 
vaccination, sanitation, and treatment as control strategies, is formulated 
in \cite{Edward}. The reproduction number for the cases with single and combined 
controls is determined and compared. The authors conclude that, when one considers 
a single control measure, treatment yields the best results, followed 
by education campaigns, sanitation, and vaccination; cf. 
the numerical simulations of \cite{Edward}. Nevertheless, 
the more control strategies are considered, better results can be obtained. 
Furthermore, the authors perform a sensitivity analysis on the key parameters 
that drive the disease dynamics in order to find their relative importance 
to cholera's spread and prevalence. In \cite{Beryl}, a SIR type model 
with a class for the bacterial concentration in the environment is proposed, 
which incorporates media coverage. The existence and stability 
of the equilibria is analyzed. Numerical simulations suggest that 
the number of infections decreases faster, when media coverage is very efficient. 
So, media alert and awareness campaigns are crucial for controlling the spread of cholera.
In \cite{Sun}, a SIR type mathematical model for cholera transmission is used 
to characterize the cholera spread in China. With the purpose of avoiding cholera 
outbreaks in China, the researchers suggest to increase the immunization 
coverage rate and to make efforts for improving environmental management, 
mainly for drinking water \cite{Sun}. 
Another important topic to have a full overview 
of mathematical modeling, within the scope of SIR models in biomathematics, 
concerns uncertainty quantification using randomized models, 
as, e.g., in \cite{Zhou, Calatayud}.
With respect to cholera, one can find, e.g., the paper \cite{Liu}, 
where the authors propose and analyze 
a stochastic epidemic model to study such disease.

The worst cholera outbreak in human history has occurred in Yemen 
and it has provoked 1 115 378 suspected cholera cases 
and 2 310 deaths, from 27 April 2017 to 1 July 2018 \cite{WHO_Yemen}.
In \cite{Nishiura:2017}, the authors forecast this cholera epidemic, 
explicitly addressing the reporting delay and ascertainment bias (see also \cite{Nishiura:2018}).
The Yemen outbreak data available in the website of World Health Organization 
is also fitted by He et al. in \cite{He}, considering a mathematical model 
based on differential equations. Their model considers five classes: 
$S$ (susceptible individuals), $I$ (infectious individuals), $R$ (recovered individuals), 
$B$ (concentration of bacterium in the environment), and $M$ (availability of medical 
resources in the country). Such model translates the interaction between the human hosts 
and the pathogenic bacteria, under the impact of limited medical resources. The results 
obtained in \cite{He} support that improvement of the public health system 
and strategic implementation of control measures with respect to time and location 
can facilitate the prevention and intervention related to this disease in Yemen.
In \cite{Lemos-Paiao2}, a SITRV (Susceptible--Infectious--Treated--Vaccinated--Recovered) 
type model with a class for bacterial concentration is proposed and mathematically studied. 
The cholera outbreak in Yemen is numerically simulated through a sub-model of the model 
mentioned previously. Moreover, a numerical simulation, which considers an hypothetical 
situation with vaccination since the beginning of the epidemic, is carried out. 
The obtained results support the importance of vaccination to prevent cholera spread.

Cholera is one of the international quarantine infectious diseases, 
as stipulated by the International Health Regulations (IHR) \cite{Cui}.
In agreement with this fact, in \cite{Lemos-Paiao} a SIQR 
(Susceptible--Infectious--Quarantined--Recovered)
type model is proposed and analyzed. In \cite{Lemos-Paiao} 
it is assumed that infective individuals are subject to quarantine 
during the treatment period. As the symptoms of the disease can not appear 
immediately after the infection (see \cite{TimeDelayCholera}),
we propose here an improvement of the work done in \cite{Lemos-Paiao},
by considering and analyzing a delayed SIQR type model. More precisely,
we introduce a discrete time-delay that represents the time between 
the instant at which an individual becomes infected and the instant 
at which he begins to have symptoms of cholera.

Several cholera outbreaks have occurred since 2007,
namely in Angola, Haiti, Zimbabwe, and Yemen \cite{Shuai,AsianScientist}.
In Haiti, the first cases of cholera happened in the Department
of Artibonite on 14th October 2010. The disease spread along the Artibonite river
and reached several departments. Only within one month,
all departments had reported cases in rural areas and places without
good conditions of public health \cite{WHO}.
We provide numerical simulations for the cholera outbreak in the Department
of Artibonite, from 1 November 2010 to 1 May 2011 \cite{WHO},
improving the results in \cite{Lemos-Paiao}.

We propose and analyze an optimal control problem with a state delay,
where the control function represents the fraction of infective individuals $I$ 
that will be submitted to treatment in quarantine until complete recovery. 
The objective is to find the optimal treatment strategy through quarantine 
that minimizes the number of infective individuals and the bacterial concentration,
as well as the cost of interventions associated with quarantine.
We apply the Pontryagin Minimum Principle for time-delayed optimal control problems
(see \cite{Gollmann1,Gollmann2}) with $L^1$-type cost functional. 
The delayed and non-delayed optimal control problems with $L^1$ 
cost functionals are analyzed, analytically and numerically, 
the solutions being interpreted from an epidemiological point of view.

The paper is organized as follows. In Section~\ref{Sec:model}, 
we propose a delayed model for the cholera transmission dynamics. 
In Section~\ref{Sec:mod:analysis_delay}, the delayed model is analyzed, 
proving the non-negativity of the solutions for non-negative initial 
conditions and providing the disease-free equilibrium, basic reproduction number $R_0$, 
and endemic equilibrium, when $R_0>1$. The stability of the equilibrium points 
is analyzed for positive time delays. Concretely, a stability analysis 
of the endemic equilibrium, as function of the ingestion rate of the bacteria 
through contaminated sources, is carried out. In Section~\ref{sec:num:simu}, 
we consider the cholera outbreak that occurred in the Department of Artibonite -- Haiti, 
from 1 November 2010 to 1 May 2011 \cite{WHO}, improving the numerical simulations 
done in \cite{Lemos-Paiao} by considering a positive time delay, 
treatment, and recovery. In Section~\ref{sec:ocp}, we formulate a
SIQRB (Susceptible--Infectious--Quarantined--Recovered--Bacterial) 
time-delayed optimal control problem with $L^1$ cost functional. 
Necessary optimality conditions are discussed. 
We consider different delayed control problems 
and compute \textit{extremal solutions}  using discretization 
and nonlinear programming methods.
Section~\ref{sec:conclusion} presents our main conclusions.


\section{Model formulation}
\label{Sec:model}

\begin{figure}[ht!]
\centering
\includegraphics[scale=1.0]{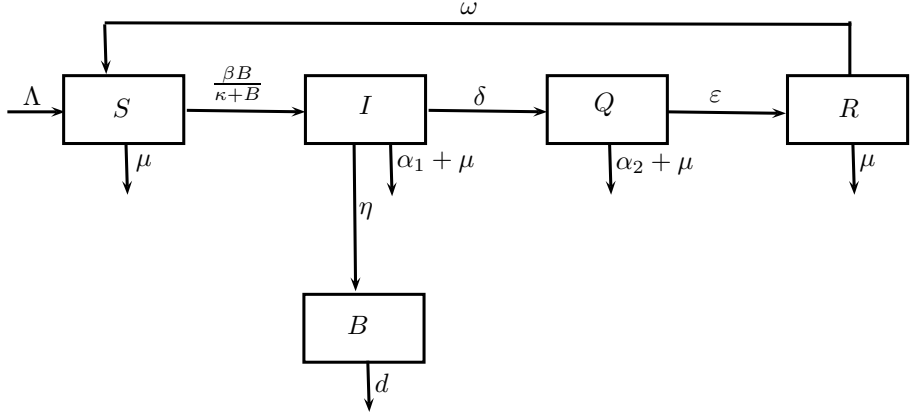}
\caption{Diagram of the dynamical model \eqref{ModeloColera_delay}.}
\label{fig:diag}
\end{figure}

We revisit the dynamic model studied in \cite{Lemos-Paiao} and introduce
a time delay $\tau\geq0$, which represents the time between the instant 
in which an individual becomes infected and the instant in which he begins 
to show symptoms. Thus, we propose a time-delayed model that involves 
a SIQR (Susceptible--Infectious--Quarantined--Recovered)
system and that also considers a class of bacterial concentration in the
dynamics of cholera. The total population $N(t)$ is divided into
four classes: susceptible $S(t)$, infectious with symptoms $I(t)$,
in treatment through quarantine $Q(t)$, and recovered $R(t)$ at time $t$,
for $t \ge 0$. Furthermore, the class $B(t)$ reflects
the bacterial concentration at time $t \ge 0$. The time delay $\tau\geq 0$ 
is related with the passage of individuals from class $S$ to class $I$. 
The introduction of this delay is done with the goal to better approximate reality.
The symptoms of cholera can appear from a few hours to 5 days after infection.
Nevertheless, they typically appear in 2--3 days \cite{TimeDelayCholera}. 
We assume that there is a positive recruitment rate $\Lambda$ into the susceptible 
class $S(t)$ and a constant natural death rate $\mu > 0$. Susceptible individuals 
can become infected with cholera at rate $\displaystyle\frac{\beta B(t)}{\kappa+B(t)}$
that is dependent on time  $t$. Note that $\beta>0$ is the ingestion
rate of the bacteria through contaminated sources, $\kappa$ is the half
saturation constant of the bacteria population, and $\displaystyle\frac{B(t)}{\kappa+B(t)}$
is the likeliness of an infective individual to have the disease with symptoms,
given a contact with contaminated sources \cite{Mwasa}. Any recovered individual
can loose the immunity at rate $\omega$ and therefore becomes susceptible again.
The infective individuals can accept to be in quarantine during a period of time.
During this time they are isolated and subject to a proper medication,
at rate $\delta$. The quarantined individuals can recover at rate
$\varepsilon$. The disease-related death rates associated with the individuals
that are infective and in quarantine are $\alpha_1$ and $\alpha_2$, respectively.
Each infective individual contributes to the increase of the bacterial concentration
at rate $\eta$. On the other hand, the bacterial concentration can decrease
at mortality rate $d$. These assumptions are translated into the following
mathematical model:
\begin{equation}
\label{ModeloColera_delay}
\begin{cases}
\begin{split}
\dot{S}(t)&=\Lambda-\displaystyle\frac{\beta B(t)}{\kappa+B(t)}S(t)
+\omega R(t)-\mu S(t),\\[0.15cm]
\dot{I}(t)&=\displaystyle\frac{\beta B(t-\tau)}{\kappa+B(t-\tau)}
S(t-\tau)-(\delta+\alpha_1+\mu)I(t),\\[0.15cm]
\dot{Q}(t)&=\delta I(t)-(\varepsilon+\alpha_2+\mu)Q(t),\\[0.15cm]
\dot{R}(t)&=\varepsilon Q(t)-(\omega+\mu)R(t),\\[0.15cm]
\dot{B}(t)&=\eta I(t)-dB(t).\\[0.1cm]
\end{split}
\end{cases}
\end{equation}
In Figure~\ref{fig:diag}, the model \eqref{ModeloColera_delay} 
is presented in a schematic way.

\begin{rmrk}
Note that 
\begin{enumerate}[i.]
\item the parameters $\Lambda$, $\mu$, $\beta$, $\kappa$, $\eta$, and $d$ are strictly positive 
while $\omega$, $\delta$, $\varepsilon$, $\alpha_1$, $\alpha_2$, and $\tau$ are nonnegative. 
	
\item individuals in class $R$ have already finalized their treatment, 
being recovered and, consequently, immune to cholera.
Following \cite{Cui,Shuai,Mwasa,Beryl}, we assume that there is no cholera induced-deaths 
for recovered individuals. In this way, it makes sense to consider that both susceptible 
and recovered individuals die at natural death rate $\mu$. In other words, we are considering 
that only non-treated infectious individuals with symptoms, 
or in treatment, can die due to cholera.
\end{enumerate}
\end{rmrk}

\section{Model analysis}
\label{Sec:mod:analysis_delay}

First, we prove that the model \eqref{ModeloColera_delay} 
makes sense from a biological point of view, because the 
solutions of \eqref{ModeloColera_delay} are non-negative 
under non-negative initial conditions. Secondly, we give 
the expressions of the disease-free and endemic equilibrium points, 
as well as the expression of the basic reproduction number.
Then, we proceed with the linearization of the model, which 
allows to derive some important results needed in the stability study  
of the equilibria. In the sequel, we shall consider the following notations:
\begin{enumerate}
\item[1)] $a_1=\delta+\alpha_1+\mu$;
\item[2)] $a_2=\varepsilon+\alpha_2+\mu$;
\item[3)] $a_3=\omega+\mu$;
\item[4)] $\rho=\Lambda\eta a_2a_3+\kappa d(a_1a_2a_3-\delta\varepsilon\omega)$;
\item[5)] $\bar{D}=a_1a_2a_3\mu+\beta(a_1a_2a_3-\delta\varepsilon\omega)$;
\item[6)] $A=a_1a_2a_3$;
\item[7)] $\tilde{A}=a_1a_2a_3-\delta\varepsilon\omega$.
\end{enumerate}


\subsection{Non-negativity of solutions}
\label{subsec:positivity:solutions}

As in \cite{Lemos-Paiao}, we assume that the initial conditions
of system \eqref{ModeloColera_delay} are non-negative:
\begin{equation}
\label{eq:init:cond_delay}
S(t) = S_0 \geq 0 \, , \quad I(t) = I_0 \geq 0 \, ,
\quad Q(t)= Q_0 \geq 0 \, , \quad R(t) = R_0 \geq 0 \, 
\quad \text{and} \quad B(t) = B_0 \geq 0 \quad
\text{ for all }\quad t\in[-\tau,0].
\end{equation}
The model studied in \cite{Lemos-Paiao} has biological meaning.
Lemma~\ref{lemma:01} proves that the corresponding delayed 
model \eqref{ModeloColera_delay}--\eqref{eq:init:cond_delay} 
also makes sense from the biological point of view.

\begin{rmrk}
One could consider suitable generic continuous functions 
for initial conditions \eqref{eq:init:cond_delay}.
However, for the concrete situation of the Haiti outbreak under investigation, 
there is no way to know such continuous functions. For this reason, 
those continuous functions are here approximated by the
constant functions \eqref{eq:init:cond_delay}.
\end{rmrk}

\begin{lmm}
\label{lemma:01}
The solutions $\left(S(t), I(t), Q(t), R(t), B(t)\right)$ 
of \eqref{ModeloColera_delay} are non-negative 
for all $t \geq -\tau$ with non-negative initial conditions
\eqref{eq:init:cond_delay}.
\end{lmm}

\begin{proof}
We have
\begin{equation*}
\begin{cases}
\begin{split}
\frac{d S(t)}{d t}\bigg|_{\xi(S)}
&=\ \Lambda + \omega R(t) > 0 \, ,\\[0.15cm]
\frac{d I(t)}{d t}\bigg|_{\xi(I)}
&=\ \displaystyle\frac{\beta B(t-\tau)S(t-\tau)}{\kappa+B(t-\tau)}
\ \geq 0 \, ,\\[0.15cm]
\frac{d Q(t)}{d t}\bigg|_{\xi(Q)}
&=\ \delta I(t)\ \geq 0 \, , \\[0.15cm]
\frac{d R(t)}{d t}\bigg|_{\xi(R)}
&=\ \varepsilon Q(t)\ \geq 0 \, , \\[0.15cm]
\frac{d B(t)}{d t}\bigg|_{\xi(B)}
&=\ \eta I(t)\ \geq 0 \, ,
\end{split}
\end{cases}
\end{equation*}
where $\xi(\upsilon)=\big \{\upsilon(t)=0 \text{ and }
S, I, Q, R, B \in C\big([-\tau,+\infty[,\mathbb{R}_0^+\big )\big \}$
and $\upsilon \in \{S, I, Q, R, B\}$.
Therefore, due to Lemma~2 in \cite{Yang:CMA:1996},
any solution of system \eqref{ModeloColera_delay} is such that
$\big(S(t), I(t), Q(t), R(t), B(t)\big ) \in \big (\mathbb{R}_0^+\big )^5$
for all $t \geq -\tau$.
\end{proof}


\subsection{Equilibrium points and the basic reproduction number} 
\label{subsec:equilibria}

When a system of ordinary differential equations with constant parameters 
is in equilibrium, their states are represented by constant functions.
In our case, in a state of equilibrium we have 
$$
\dot{S}(t) = \dot{I}(t) = \dot{Q}(t) = \dot{R}(t) = \dot{B}(t) = 0.
$$
Consequently, we get constant functions. 
Therefore, time delays do not affect the expressions of equilibria.
Thus, equilibria of delayed and non-delayed $SIQRB$ models for cholera, 
with constant parameters, are the same. 
From \cite{Lemos-Paiao}, we know that model 
\eqref{ModeloColera_delay} has a disease-free equilibrium (DFE) given by
\begin{equation}
\label{eq:DFE}
E^0=\big(S^0,I^0,Q^0,R^0,B^0\big)
=\left(\frac{\Lambda}{\mu},0,0,0,0\right)
\end{equation}
and that the basic reproduction number has the following expression:
\begin{equation}
\label{R0}
R_0=\frac{\beta\Lambda\eta}{\mu\kappa da_1}.
\end{equation}
Moreover, when $R_0>1$, there is an endemic equilibrium given by
\begin{equation}
\label{EndemicEquilibrium}
E^*=(S^*,I^*,Q^*,R^*,B^*),
\end{equation}
where
\begin{equation*}
S^*=\displaystyle\frac{a_1\rho}{\eta\bar{D}}, \quad
I^*=\displaystyle\frac{\beta\Lambda a_2a_3(R_0-1)}{R_0\bar{D}}, \quad
Q^*=\displaystyle\frac{\beta\Lambda a_3\delta(R_0-1)}{R_0\bar{D}}, \quad
R^*=\displaystyle\frac{\beta\Lambda\delta\varepsilon(R_0-1)}{R_0\bar{D}}, 
\quad B^*=\displaystyle\frac{\beta\Lambda\eta a_2a_3(R_0-1)}{R_0\bar{D}d}.
\end{equation*}


\subsection{Linearisation of the model}
\label{subsec:linearisation:model}

Let us denote the state variables by $x_1=S$, $x_2=I$, $x_3=Q$, $x_4=R$ and $x_5=B$.
Then, we can write the system \eqref{ModeloColera_delay} in the following way:
\begin{equation*}
\dot{x}(t)=f\big (x(t),x(t-\tau)\big),
\end{equation*}
where $x=\big(x_1,x_2,x_3,x_4,x_5\big)$. 
Let $\bar{E}=(\bar{x}_1,\bar{x}_2,\bar{x}_3,\bar{x}_4,\bar{x}_5)$ 
be an arbitrary equilibrium point of \eqref{ModeloColera_delay} 
and let us consider the following change of variables:
\begin{equation*}
z_i(t)=x_i(t)-\bar{x}_i, \quad i=1,\dots,5.
\end{equation*}
Thus, the linearized system of \eqref{ModeloColera_delay} is given by
\begin{equation*}
\dot{z}=\frac{\partial f}{\partial x}\bigg|_{\bar{E}}z
+\frac{\partial f}{\partial x_{\tau}}\bigg|_{\bar{E}}z_{\tau},
\end{equation*}
where $z(t)=\big(z_1(t),z_2(t),z_3(t),z_4(t),z_5(t)\big)$, 
$z_{\tau}(t)=z(t-\tau)$ and $x_{\tau}(t)=x(t-\tau)$. Furthermore, we have
\begin{equation*}
A_0:=\frac{\partial f}{\partial x}\bigg|_{\bar{E}}
=\left[
\begin{matrix}
-\bar{\lambda}-\mu & 0 & 0 & \omega & -C\\
0 & -a_1 & 0 & 0 & 0\\
0 & \delta & -a_2 & 0 & 0\\
0 & 0 & \varepsilon & -a_3 & 0\\
0 & \eta & 0 & 0 & -d
\end{matrix}
\right]
\end{equation*}
and
\begin{equation*}
A_1:=\frac{\partial f}{\partial x_{\tau}}\bigg|_{\bar{E}}
=\left[
\begin{matrix}
0 & 0 & 0 & 0 & 0\\
\bar{\lambda} & 0 & 0 & 0 & C\\
0 & 0 & 0 & 0 & 0\\
0 & 0 & 0 & 0 & 0\\
0 & 0 & 0 & 0 & 0
\end{matrix}
\right],
\end{equation*}
where $\bar{\lambda}=\displaystyle\frac{\beta\bar{x}_5}{\kappa+\bar{x}_5}$ 
and $C=\displaystyle\frac{\beta\kappa\bar{x}_1}{(\kappa+\bar{x}_5)^2}$. So, 
the characteristic polynomial associated with the linearized system of model 
\eqref{ModeloColera_delay} is given by
\begin{equation}
\label{caracteristic_polynomial_delay}
p_{\tau}(\chi)=\det(\chi I_5-A_0-e^{-\tau\chi}A_1)
=P_1(\chi)+e^{-\tau \chi}P_2(\chi),
\end{equation}
where 
$$
P_1(\chi)=(\chi+a_1)(\chi+a_2)(\chi+a_3)(\chi+d)(\chi+\bar{\lambda}+\mu)
$$
and
$$
P_2(\chi)=-\eta C\chi^3-\eta C(a_2+a_3+\mu)\chi^2
-\Big(\eta C(a_2a_3+a_2\mu+a_3\mu)+\delta\varepsilon\omega\bar{\lambda}\Big)\chi
-\eta Ca_2a_3\mu-\delta\varepsilon\omega d\bar{\lambda}.
$$


\subsection{Stability analysis}
\label{subsec:stability}

In order to study the stability of the equilibria, we are going 
to follow the approach of \cite{Cooke}, corrected in \cite{Boese} 
(see also Theorem~4.1 of \cite[p.~83]{Kuang}). We shall prove that
the conditions $(i)$, $(ii)$, $(iv)$ and $(v)$ of Theorem~4.1 
of \cite[p.~83]{Kuang} are satisfied for an arbitrary equilibrium point 
$$
\bar{E}=(\bar{x}_1,\bar{x}_2,\bar{x}_3,\bar{x}_4,\bar{x}_5)\in(\mathbb{R}_0^+)^5
$$ 
of model \eqref{ModeloColera_delay}. The polynomial $P_1$ has only real zeros 
($\chi=-a_1$ or $\chi=-a_2$ or $\chi=-a_3$ or $\chi=-d$ or $\chi=-\bar{\lambda}-\mu$).
Thus, the polynomials $P_1$ and $P_2$ can not have common imaginary zeros 
and condition $(i)$ is satisfied. In order to
fulfill the hypothesis of condition $(ii)$, we are going 
to compute $P_1(yi)$ and $P_2(yi)$. We have
\begin{equation*}
\begin{split}
P_1(yi)&=\ \big [a_1(\bar{\lambda}+\mu)-y^2\big ]
\big [a_2a_3d-(a_2+a_3+d)y^2\big ]
-(a_1+\bar{\lambda}+\mu)\big [a_2a_3+(a_2+a_3)d-y^2\big ]y^2\\
&\quad +\Big \{\big [a_1(\bar{\lambda}+\mu)-y^2\big ]
\big[a_2a_3+(a_2+a_3)d-y^2\big ]+(a_1+\bar{\lambda}+\mu)
\big[a_2a_3d-(a_2+a_3+d)y^2\big ]\Big \}yi\\
&=\ \overline{P_1(-yi)}
\end{split}
\end{equation*}
and one obtains that
\begin{equation*}
\begin{split}
P_2(yi)&=\eta C(a_2+a_3+\mu)y^2-\eta Ca_2a_3\mu
-\delta\varepsilon\omega\bar{\lambda}d+\Big \{\eta Cy^3
-\big [\eta C(a_2a_3+a_2\mu+a_3\mu)
+\delta\varepsilon\omega\bar{\lambda}\big ]y\Big \}i\\
&=\ \overline{P_2(-yi)}.
\end{split}
\end{equation*}
Therefore, condition $(ii)$ is also satisfied. 
As the degree of polynomial $P_1$, equal to five,
is bigger than the degree of $P_2$, equal to three, 
then the condition $(iv)$, given by
$$
\lim\limits_{|\lambda|\rightarrow\infty,\ \Re(\lambda)\geq0} 
\sup \left\{\left|\frac{P_2(\lambda)}{P_1(\lambda)}\right|\right\}<1,
$$
is obviously satisfied. Furthermore, the function defined by 
$F(y)=|P_1(yi)|^2-|P_2(yi)|^2$ is a polynomial with degree equal to ten. 
Thus, the function $F$ has at most a finite number (ten) of real zeros. 
Concluding, the condition $(v)$ is also verified.


\subsubsection{Disease-free equilibrium}

Now, we are going to study the stability of the 
DFE \eqref{eq:DFE} of the delayed model \eqref{ModeloColera_delay}.

\begin{thrm}[Stability of the DFE \eqref{eq:DFE}]
\label{thm:01}
Assume that $R_0\neq1$. If $a_1d<1$, then
there exists $\tau^*\in\mathbb{R}_0^+$ such that:
\begin{itemize}
\item there is at most a finite number of stability switches, 
when $\tau\in[0,\tau^*]$;
		
\item instability occurs, when $\tau \in\, ]\tau^*,+\infty[$.
\end{itemize}
For all $\tau\geq0$, if $a_1d\geq1$, then the DFE \eqref{eq:DFE} is:
\begin{itemize}
\item locally asymptotic stable, when $R_0<1$;
\item unstable, when $R_0>1$.
\end{itemize}
\end{thrm}

\begin{proof}
In order to study the stability of the DFE, we follow the approach of \cite{Cooke,Boese}. 
We already know that the conditions $(i)$, $(ii)$, $(iv)$, and $(v)$ of 
Theorem~4.1 of \cite[p.~83]{Kuang} are satisfied for the DFE. Thus, we  
analyze now the condition $(iii)$ and compute the zeros of the polynomial $F$ for the DFE. 
Computing $P_1(0)+P_2(0)$, we obtain
\begin{equation*}
P_1(0)+P_2(0)=a_1a_2a_3d\mu-\eta Ca_2a_3\mu
= a_1a_2a_3d\mu-a_1dR_0a_2a_3\mu
= a_1a_2a_3d\mu(1-R_0).
\end{equation*}
Concluding, $P_1(0)+P_2(0)\neq 0 \Leftrightarrow R_0\neq1$, 
because $a_1,a_2,a_3,d,\mu>0$. Therefore, the condition $(iii)$ 
is verified for the DFE if and only if $R_0\neq1$.
	
As the conditions $(i)$--$(v)$ are verified with respect to the DFE 
if and only if $R_0\neq1$ holds, the stability of the  DFE depends 
on the roots of the polynomial $F$, according to Theorem 4.1 of \cite[p.~83]{Kuang}. 
Solving $F(y)=0$, we get
\begin{gather*}
y=\pm a_2i \vee y=\pm a_3i \vee y=\pm\mu i \vee y
=\pm\frac{\sqrt{2}}{2}\left(\sqrt{a_1^2+d^2+\sqrt{(a_1^2-d^2)^2+4}}\right)i\\
\vee y=\pm\frac{\sqrt{2}}{2}\sqrt{-a_1^2-d^2+\sqrt{(a_1^2-d^2)^2+4}}.
\end{gather*}
If $-a_1^2-d^2+\sqrt{(a_1^2-d^2)^2+4}>0$, which is equivalent to
\begin{align*}
&\sqrt{(a_1^2-d^2)^2+4}>a_1^2+d^2\underset{a_1^2+d^2>0}{\Leftrightarrow} 
(a_1^2-d^2)^2+4 > (a_1^2+d^2)^2\Leftrightarrow (a_1d)^2
<1\underset{a_1,d>0}{\Leftrightarrow} 0<a_1d<1,
\end{align*}
then the polynomial $F$ has at least one positive root, which is simple. 
According to Theorem 4.1 of \cite[p.~83]{Kuang}, when $a_1d<1$, 
we can state that there is $\tau^*>0$ such that
\begin{itemize}
\item at most a finite number of stability switches may occur, 
if $\tau\in[0,\tau^*]$;
\item instability occurs, if $\tau \in \, ]\tau^*,+\infty[$.
\end{itemize}
On the other hand, if $a_1d\geq1$, then the polynomial 
$F$ has no positive roots. In this case, according with 
item (a) of Theorem~4.1 of \cite[p.~83]{Kuang},
the stability/instability is determined by the 
corresponding stability/instability that occurs when $\tau=0$. 
When $\tau=0$, one has the model studied in \cite{Lemos-Paiao}. 
Therefore, when $a_1d\geq1$, the DFE is
locally asymptotic stable, if $R_0<1$;
unstable, if $R_0>1$; for all $\tau\geq0$.
\end{proof}

\begin{rmrk}
The computation of an analytic expression for $\tau^*$, 
whose existence is proved in Theorem~\ref{thm:01}, is cumbersome.
Our problem is complex and it is very difficult, or even impossible, 
to obtain some analytical results. We offer as an open problem 
the question of how to compute the exact expression of $\tau^*$.
Eventually, the procedure of reference \cite{Kuang} can be helpful.
\end{rmrk}


\subsubsection{Endemic equilibrium}

For the endemic equilibrium point \eqref{EndemicEquilibrium}, we have
\begin{equation*}
C^*=\frac{\beta\kappa S^*}{(\kappa+B^*)^2}
=\frac{\beta\Lambda\Big (R_0kd(a_1a_2a_3-\delta\varepsilon\omega)
+\Lambda\eta a_2a_3\Big )}{\mu kR_0^2\rho},
\end{equation*}
which implies
\begin{equation*}
\begin{split}
\eta C^*&=\frac{\mu kda_1R_0\Big (R_0kd(a_1a_2a_3
-\delta\varepsilon\omega)+\Lambda\eta a_2a_3\Big)}{\mu 
kR_0^2\rho}=\frac{a_1d}{R_0\rho}\Big (R_0kd(a_1a_2a_3
-\delta\varepsilon\omega)+\Lambda\eta a_2a_3\Big )\\
&=\frac{a_1d}{R_0\rho}\Big (R_0(\rho-\Lambda\eta a_2a_3)+\Lambda\eta a_2a_3\Big )
=\frac{a_1d}{R_0\rho}\Big (R_0\rho+\Lambda\eta a_2a_3(1-R_0)\Big )\\
&=a_1d+\frac{\Lambda\eta a_1a_2a_3d(1-R_0)}{R_0\rho}
=a_1d+\frac{\Lambda\eta a_1^2a_2a_3d^2\mu\kappa(1-R_0)}{\beta\Lambda\eta\rho}\\
&=a_1d+\frac{(a_1d)^2a_2a_3\mu\kappa(1-R_0)}{\beta\rho}
=a_1d\left(1+\frac{a_1a_2a_3\mu\kappa d(1-R_0)}{\beta\rho}\right).
\end{split}
\end{equation*}
Again, we have to study the condition $(iii)$ of 
Theorem 4.1 of \cite[p.~83]{Kuang} and the roots of $F$ 
with respect to the endemic equilibrium \eqref{EndemicEquilibrium}. 
Computing $P_1(0)+P_2(0)$ for $E^*$, we obtain
\begin{equation*}
\begin{split}
P_1(0)+P_2(0)
&=-\frac{a_2a_3}{kR_0^2\rho}\big (\beta\Lambda\eta-R_0^2\mu kda_1\big)
\big (R_0kd(a_1a_2a_3-\delta\varepsilon\omega)+\Lambda\eta a_2a_3\big )\\
&=-\frac{a_2a_3}{kR_0^2\rho}\big (R_0\mu kda_1-R_0^2\mu kda_1\big)
\big(R_0kd(a_1a_2a_3-\delta\varepsilon\omega)+\Lambda\eta a_2a_3\big )\\
&=-\frac{R_0\mu kda_1a_2a_3}{kR_0^2\rho}\big (1-R_0\big)
\big(R_0kd(a_1a_2a_3-\delta\varepsilon\omega)+\Lambda\eta a_2a_3\big )\\
&=-\frac{\mu da_1a_2a_3}{R_0\rho}\big (1-R_0\big )
\big(R_0kd(a_1a_2a_3-\delta\varepsilon\omega)+\Lambda\eta a_2a_3\big).
\end{split}
\end{equation*}
As $\mu da_1a_2a_3>0$, $R_0\rho>0$, and 
$R_0kd(a_1a_2a_3-\delta\varepsilon\omega)+\Lambda\eta a_2a_3>0$, 
then $P_1(0)+P_2(0)\neq0$ with respect to $E^*$ if and only if $R_0\neq1$. 
Now, we are going to write function $F$ as a function of $E^*$  
in order to obtain its roots. We obtain
\begin{equation}
\label{eq:Fend}
F_{end}(y)=F(y)|_{E^*}=c_0+c_2y^2+c_4y^4+c_6y^6+c_8y^8+c_{10}y^{10},
\end{equation}
where $c_0$, $c_2$, $c_4$, $c_6$, $c_8$, and $c_{10}$ 
are real coefficients given by
\begin{eqnarray*}
c_0&=&\ \left(\frac{A\mu}{\beta\rho}\right)^2d^3k(A\mu
+\beta\tilde{A})(R_0-1)\left\{\frac{\beta^2\Lambda\eta}{\mu}\left(a_2
a_3+\frac{\delta\varepsilon\omega}{a_1}\right)+\rho
+kd(A\mu-2\beta\delta\varepsilon\omega)\right\};\\[2ex]
c_2&=&\ \left\{(a_2a_3\mu)^2+\frac{2a_1(a_2a_3)^3
\mu^2kd(R_0-1)}{\rho}\right\}(a_1^2+d^2)\\
&& +\left\{\frac{(a_1d)^2a_2a_3\mu k(R_0-1)}{\beta\rho}\right\}
\left\{\frac{a_2a_3\mu k(R_0-1)}{\beta\rho}\Big [A^2(\beta^2-d^2)
+(a_1d)^2(\beta^2-\mu^2)(a_2^2+a_3^2)\right.\\
&& \left.+2\beta\delta\varepsilon\omega(a_2a_3+a_2\mu
+a_3\mu-a_2d-a_3d-\mu d)a_1d+\beta^2\big [(a_2a_3d)^2
-(\delta\varepsilon\omega)^2\big ]\Big ]\right.\\
&& -2\beta\delta\varepsilon\omega(a_2a_3+a_2\mu+a_3\mu-a_2d-a_3d-\mu d)\\
&& +2a_1d\Big [(a_2a_3)^2+(a_2\mu)^2+(a_3\mu)^2+\beta\mu(a_2^2+a_3^2)\Big ]\bigg\};\\
c_4&=&\ \frac{(a_1d)^2a_2a_3\mu\kappa(R_0-1)}{\beta\rho}\bigg\{
\frac{a_1d(a_2^2+a_3^2+\mu^2)}{\beta\rho}(\beta\rho 
+ \beta\kappa d\tilde{A} + A\mu\kappa d) + 2a_1\mu d\beta\\
&& + \frac{2\kappa d\delta\varepsilon\omega(A\mu+\beta \tilde{A})}{\rho}\bigg\} 
+ (a_2a_3)^2(a_1^2+d^2) + (a_1d\lambda^*)^2\\
&& + \{(a_1a_2)^2 + (a_1a_3)^2 + (a_2a_3)^2 
+ (a_2d)^2 + (a_3d)^2\}(\lambda^*+\mu)^2;\\[2ex]
c_6&=&\ \frac{(a_1d)^3a_2a_3\mu\kappa(R_0-1)}{(\beta\rho)^2}(
\beta\rho + \beta\kappa d\tilde{A}+A\mu\kappa d)
+a_2^2(a_1^2+d^2)+a_3^2(a_1^2+a_2^2+d^2)\\
&& +(\lambda^*+\mu)^2(a_1^2+a_2^2+a_3^2+d^2);\\[2ex]
c_8&=&\ a_1^2+a_2^2+a_3^2+d^2+(\lambda^*+\mu)^2;\\[2ex]
c_{10}&=&1;
\end{eqnarray*}
with $\lambda^*=\displaystyle\frac{a_1a_2a_3\mu\kappa d(R_0-1)}{\rho}$. 
The study of the stability of $E^*$ depends on the roots of the polynomial 
$F_{end}$ \eqref{eq:Fend}. It is not easy to obtain their analytical expressions, 
but we can note that if $R_0>1$, then the coefficients $c_4$, $c_6$, $c_8$, 
and $c_{10}$ are all non-negative for any admissible parameters. Nevertheless, 
the sign of the coefficients $c_0$ and $c_2$ is yet an open question. 
If all coefficients would be positive, then the polynomial $F_{end}$ 
would not have positive roots by \emph{Descartes' Rule of Signs}. Thus, 
the stability/instability would be determined by the stability/instability 
that occurs when $\tau=0$, according to item (a) of 
Theorem 4.1 of \cite[p.~83]{Kuang}. In this way, we would obtain 
the stability result expressed in Theorem~6 of \cite{Lemos-Paiao} 
for the endemic equilibrium \eqref{EndemicEquilibrium} 
of the delayed model \eqref{ModeloColera_delay}.
Though $c_0$ and $c_2$ are given by complicated expressions, 
we can derive some conclusions about their signs by studying 
them as a function of the ingestion rate $\beta\in\ ]0,5]$ 
and fixing all parameters to the values of Table~\ref{Tab_Parameter} 
that is introduced in Section~\ref{sec:num:simu} of numerical 
simulations associated with SIQRB delayed model \eqref{ModeloColera_delay}.
Thus, for the existence of a disease we have to assume that $R_0>1$, 
which is equivalent to $\beta>4.772690\times10^{-2}$.
Analyzing the signs of the coefficients $c_0$ and $c_2$ 
as functions of $\beta\in\ ]0,5]$, we obtain
\begin{equation*}
\begin{cases}
c_0=0\Leftrightarrow\beta\simeq 2.698643\times10^{-5} 
\vee \beta\simeq 2.468318\times10^{-2} 
\vee \beta\simeq 4.772690\times10^{-2}; \\[0.1cm]
c_2=0 \Leftrightarrow \beta\simeq 4.772655\times10^{-2}.
\end{cases}
\end{equation*}
Numerically, we analyze 
the signs of the coefficients $c_0$ and $c_2$ 
for $\beta\in\ ]0,5]$. We conclude that
\begin{equation*}
\begin{cases}
c_0>0 \Leftrightarrow \beta\in\ ]2.698643\times10^{-5},2.468318\times10^{-2}[\ 
\cup\ ]4.772690\times10^{-2},5];\\[0.1cm]
c_2>0 \Leftrightarrow \beta\in\ ]4.772655\times10^{-2},5].
\end{cases}
\end{equation*}
Let us consider $4.8\times10^{-2}\leq\beta\leq5$ 
and the other parameters fixed to the values of Table~\ref{Tab_Parameter}. 
Hence, a disease occurs in view of $R_0>1$. In this case, $F_{end}$ has only 
positive coefficients. So, by \emph{Descartes' Rule of Signs}, we can state 
that the polynomial $F_{end}$ has no positive roots. Therefore, using the 
Theorem 4.1 of \cite[p.~83]{Kuang}, we can conclude that the stability/instability 
for $\tau=0$ remains for all $\tau\geq0$. As $R_0>1$, we can state that, 
for $\beta\in[4.8\times10^{-2},5]$, the endemic equilibrium obtained with 
the values of Table~\ref{Tab_Parameter} is locally asymptotic stable 
(see Theorem~6 of \cite{Lemos-Paiao}).


\section{Numerical simulations of the SIQRB delayed model}
\label{sec:num:simu}

The cholera outbreak that occurred in the Department of Artibonite -- Haiti, 
from $1^{\text{st}}$ November 2010 to $1^{\text{st}}$ May 2011 (see \cite{WHO}), 
is approximated in \cite{Lemos-Paiao} by numerical simulations of model 
\eqref{ModeloColera_delay} in the particular case when
$\tau=\omega=\delta=\varepsilon=\alpha_2=Q(0)=R(0)=0$,
which originates the sub-model given by
\begin{equation}
\label{SubModeloColera}
\begin{cases}
\begin{split}
\dot{S}(t)&=\Lambda-\displaystyle\frac{\beta B(t)}{\kappa+B(t)}S(t)-\mu S(t),\\[0.15cm]
\dot{I}(t)&=\displaystyle\frac{\beta B(t)}{\kappa+B(t)}S(t)-(\alpha_1+\mu)I(t),\\[0.15cm]
\dot{B}(t)&=\eta I(t)-dB(t).
\end{split}
\end{cases}
\end{equation}
Here we improve the results of \cite{Lemos-Paiao}
by considering system \eqref{ModeloColera_delay} instead of \eqref{SubModeloColera}.
As there are some parameters for which we do not find 
the adequate values for them in the literature, we now fit the real cholera data 
of Haiti by considering treatment, recovery, and time delay through the full SIQRB 
delayed model \eqref{ModeloColera_delay}, with the purpose 
to approximate the parameter values of $\tau, \delta, \beta$ and $\alpha_1$.

We use a best-fitting procedure that considers 
the full SIQRB delayed model. It consists in searching the values of parameters
\begin{enumerate}
\item[i)] $\tau \in [\tau_{\min},\tau_{\max}]$;
\item[ii)] $\delta \in [\delta_{\min},\delta_{\max}]$;
\item[iii)] $\beta \in [\beta_{\min},\beta_{\max}]$;
\item[iv)] $\alpha_1 \in [\alpha_{1_{\min}},\alpha_{1_{\min}}]$;
\end{enumerate}
while minimizing the quadratic error between the real data and the number 
of infected individuals predicted by the SIQRB delayed model. 
We have considered
\begin{enumerate}
\item[i)] $\tau_{\min} = 2$ and $\tau_{\max}=3$, since symptoms 
typically appear in 2--3 days (see \cite{TimeDelayCholera}); 
\item[ii)] $\delta_{\min} = 0.01$ and $\delta_{\max} = 0.02$; 
\item[iii)] $\beta_{\min} = 0.7$ and $\beta_{\max} = 1.2$; 
\item[iv)] $\alpha_{1_{\min}} = 0.005$ and $\alpha_{1_{\max}} = 0.025$. 
\end{enumerate}
\begin{figure}[ht!]
\centering
\includegraphics[scale=0.4]{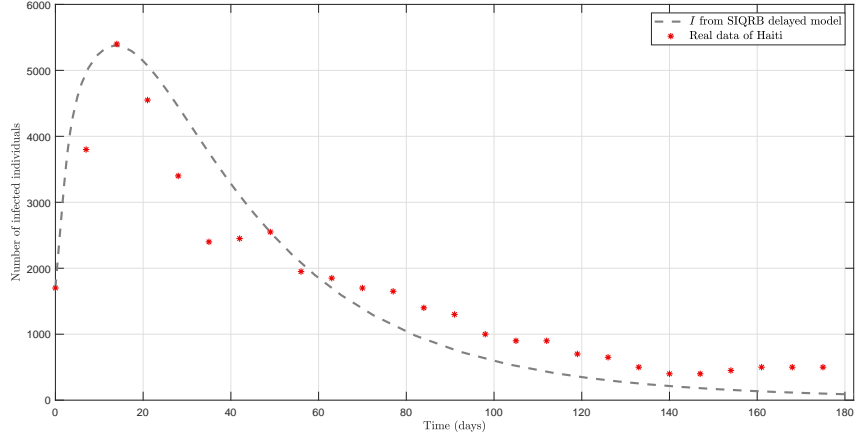}
\caption{Numerical solution $I$ through the SIQRB delayed model 
for values of Table~\ref{Tab_Parameter} with $\tau = 2.0$, 
$\delta = 0.020$, $\beta = 0.70$ and $\alpha_1 = 0.0120$
(dashed grey line) versus real data of Haiti (red stars).}
\label{FigHaiti_semQ_undelayed_delayed}
\end{figure}
The best fitted values found were: $\tau = 2.0$, $\delta = 0.020$, 
$\beta = 0.70$ and  $\alpha_1 = 0.0120$ (see Figure~\ref{FigHaiti_semQ_undelayed_delayed}).  
Note that all the other considered parameter values, with respect to the SIQRB delayed model, 
are presented in Table~\ref{Tab_Parameter}. All the numerical computations were carried out 
in \textsf{Matlab} with the help of its routines for the integration of model 
\eqref{ModeloColera_delay}.
\begin{table}[ht!]
\doublerulesep 0.1pt
\tabcolsep 7.8mm
\centering
\vspace*{2mm}
\renewcommand{\arraystretch}{1.3}
\setlength{\tabcolsep}{10pt}
\footnotesize
\caption{\rm Parameter values and initial conditions for the time-delayed 
SIQRB control model \eqref{ModeloColeraControlo}.}
\label{Tab_Parameter}
{\begin{tabular*}{15.5cm}{cccc}
\hline\hline\hline
\raisebox{-2ex}[0pt][0pt]{Parameter} & \raisebox{-2ex}[0pt][0pt]{Description} 
& \raisebox{-2ex}[0pt][0pt]{Value} & \raisebox{-2ex}[0pt][0pt]{Reference} \\ \\ \hline
\small $\Lambda$ & \small Recruitment rate & \small $24.4N(0)/365\ 000$ (day$^{-1}$) 
& \small \cite{BirthRate}\\
\small $\mu$ & \small Natural death rate & \small 2.2493$\times10^{-5}$ (day$^{-1}$) 
& \small \cite{DeathRate}\\
\small $\beta$ & \small Ingestion rate & \small 0.7 (day$^{-1}$)  & \small Best fitted\\
\small $\kappa$ & \small Half saturation constant & \small$10^6$ (cell/ml) & \small \cite{Sanches}\\
\small $\omega$ & \small Immunity waning rate & \small 0.4/365 (day$^{-1}$) & \small \cite{Neilan}\\
\small $\delta$ & \small Quarantine rate & \small 0.02 (day$^{-1}$) & \small Best fitted\\
\small $\varepsilon$ & \small Recovery rate & \small 0.2 (day$^{-1}$) & \small \cite{Mwasa}\\
\small $\alpha_1$ & \small Death rate (infective) & \small 0.012 (day$^{-1}$) & \small Best fitted\\
\small $\alpha_2$ & \small Death rate (quarantined)& \small 0.0001 (day$^{-1}$) & \small \cite{Mwasa}\\
\small $\eta$ & \small Shedding rate (infective) & \small 10 (cell/ml day$^{-1}$ person$^{-1}$) 
& \small \cite{Capone}\\
\small $d$ & \small Bacteria death rate & \small0.33 (day$^{-1}$) & \small \cite{Capone}\\
\small $\tau$ & \small Time delay associated to the state 
& \small 2 (days) & \small Best fitted in agreement with \cite{TimeDelayCholera}\\
\small $S(0)$ & \small Susceptible individuals at $t=0$ & \small 5750 (person) & \small Assumed \\
\small $I(0)$ & \small Infective individuals at $t=0$ & \small 1700 (person) & \small \cite{WHO}\\
\small $Q(0)$ & \small Quarantined individuals at $t=0$ & \small 0 (person) & \small Assumed \\
\small $R(0)$ & \small Recovered individuals at $t=0$ & \small 0 (person) & \small Assumed \\
\small $B(0)$ & \small Bacterial concentration at $t=0$ 
& \small$275\times10^3$ (cell/ml) & \small Assumed\\
\small $T$ & \small Final time & \small 182 (days) & \small \cite{WHO}\\
\small $W_u$ & \small Measure of the treatment cost & \small$1000$ & \small Assumed\\
\small $u_{\max}$ & \small Maximum value of the control & \small$4$ & \small Assumed\\
\hline\hline\hline
\end{tabular*}}
\renewcommand{\arraystretch}{1}
\end{table}

\section{Optimal control problem for the time-delayed SIQRB model}
\label{sec:ocp}

In the current section, we begin by formulating an optimal control 
problem related with the delayed SIQRB model \eqref{ModeloColera_delay} 
and then we apply the Pontryagin Minimum Principle. 
We end Section~\ref{sec:ocp} with the study of several numerical solutions 
of the proposed optimal control problem varying the weights 
$W_I$ and $W_B$ of the cost functional associated with classes $I$ and $B$, 
respectively. These studies are done with the purpose to obtain control 
strategies that could have stopped the spread of the outbreak 
in Haiti considered in Section~\ref{sec:num:simu}. 
In order to simplify the notation, 
we consider throughout the current section that the state variable is 
$$
X=(x_1,x_2,x_3,x_4,x_5)=(S,I,Q,R,B).
$$


\subsection{Formulation of the optimal control problem}
\label{sec:SIQRB-ocp}

We propose a time-delayed optimal control problem associated with model \eqref{ModeloColera_delay}. 
The set of differential equations of such optimal control problem is obtained by adding to model 
\eqref{ModeloColera_delay} a control function $u(\cdot)$ taking values in 
$\left[1,u_{\max}\right]$, which accelerates the movement from class $I$ to $Q$. 
In practice, by applying such control measure, 
infected individuals can be put under quarantine faster.
If $u \equiv 1$, then all infective individuals are quarantined at the end 
of $\frac{1}{\delta} = \frac{1}{0.02} = 50$ days. On the other hand, 
if $u \equiv u_{\max}$, then all infective people are put under quarantine 
at the end of $\frac{1}{\delta u_{\max}} < 50$ days. 
The model with control is given by the following system of non-linear 
and delayed ordinary differential equations:
\begin{equation}
\label{ModeloColeraControlo}
\begin{cases}
\begin{split}
\dot{S}(t)&=\Lambda-\displaystyle\frac{\beta B(t)}{\kappa
+B(t)}S(t)+\omega R(t)-\mu S(t),\\[0.15cm]
\dot{I}(t)&=\displaystyle\frac{\beta B(t-\tau)}{\kappa
+B(t-\tau)}S(t-\tau)-(\delta u(t)+\alpha_1+\mu)I(t),\\[0.15cm]
\dot{Q}(t)&=\delta u(t)I(t)-(\varepsilon+\alpha_2+\mu)Q(t),\\[0.15cm]
\dot{R}(t)&=\varepsilon Q(t)-(\omega+\mu)R(t),\\[0.15cm]
\dot{B}(t)&=\eta I(t)-dB(t),\\[0.1cm]
\end{split}
\end{cases}
\end{equation}
with initial conditions given by
\begin{equation}
\label{eq:init:cond_delay_OC}
S(t) = S_0 \geq 0 \, , \quad I(t) = I_0 \geq 0 \, ,
\quad Q(t)= Q_0 \geq 0 \, , \quad R(t) = R_0 \geq 0 \, \quad \text{and}
\quad B(t) = B_0 \geq 0 \,
\end{equation}
for all $t\in[-\tau,0]$.
Our aim is to minimize the number of infective individuals
and the bacterial concentration, as well as the cost of
interventions associated with the control treatment through quarantine. 
Thus, we consider the following objective functional:
\begin{equation}
\label{costfunction}
J\big (X(\cdot),u(\cdot)\big ) = \int_0^{T} 
\big( W_I I(t) + W_B B(t) + W_u u(t) \big ) dt \, ,
\end{equation}
where $T>0$ is the final time, $W_I>0$ and $W_B>0$ are the weights associated 
with infected individuals and bacterial concentration in the environment, 
respectively, and constant $W_u>0$ is a measure of the cost 
of interventions associated with the control $u$. It is important to emphasize 
that it is more appropriate to use a linear control integrand 
in the biomedical framework, as we do in \eqref{costfunction},
since the cost is directly proportional to the dosage control. Moreover, in real 
life, it is easier to apply a linear control measure than a quadratic one
(see \cite{Ledzewicz1,Ledzewicz2}). The set $\mathcal{X}$ 
of admissible trajectories is given by
\begin{equation*}
\mathcal{X} =
\big \{
X(\cdot) \in W^{1,1}\big ([0,T];\mathbb{R}^5\big ) \, | \,
\eqref{ModeloColeraControlo} \text{ and } 
\eqref{eq:init:cond_delay_OC} \text{ are satisfied}\big \}
\end{equation*}
and the admissible control set $\mathcal{U}$ is given by
\begin{equation*}
\mathcal{U} =
\big \{
u(\cdot) \in L^1\big ([0, T]; \mathbb{R}\big ) \, | \,
1 \leq u (t) \leq u_{\max} ,  \, \forall \, t \in [0, T] \,
\big \}.
\end{equation*}
Recall that $W^{1,1}$ and $L^1$ are, respectively, the class of absolutely 
continuous functions and the space of Lebesgue integrable functions.

The optimal control problem consists 
of determining the vector function
$$
X^\diamond(\cdot) = \big(S^\diamond(\cdot),
I^\diamond(\cdot), Q^\diamond(\cdot),
R^\diamond(\cdot), B^\diamond(\cdot)\big)
\in \mathcal{X}
$$ 
associated with an admissible control
$u^\diamond(\cdot) \in \mathcal{U}$ on the time interval $[0, T]$,
minimizing the cost functional \eqref{costfunction}, \textrm{i.e.},
\begin{equation}
\label{mincostfunct}
J\big (X^\diamond(\cdot),u^\diamond(\cdot)\big )
= \min_{(X(\cdot),u(\cdot))
\in\mathcal{X}\times\mathcal{U}} 
J\big (X(\cdot),u(\cdot)\big).
\end{equation}


\subsection{Necessary optimality conditions: the Pontryagin Minimum Principle}
\label{subsec:SIQRB-ocp:nec_cond}

The necessary optimality conditions for an optimal solution 
of \eqref{mincostfunct} are given by the Pontryagin Minimum Principle 
for time-delayed optimal control problems \cite{Gollmann1,Gollmann2}. 
Let us denote the delayed state variables in the following way:
\begin{enumerate}
\item[i)] $S(t-\tau)=x_1(t-\tau)=y_1(t)$;
\item[ii)] $B(t-\tau)=x_5(t-\tau)=y_5(t)$.
\end{enumerate}
The Hamiltonian function associated with our optimal control problem is given by
\begin{equation}
\label{Hamiltonian}
\begin{split}
H(X,y_1,y_5,u,\lambda)
&= \ W_I x_2 + W_B x_5  + W_u u+ \lambda_1 \left( \Lambda
-\displaystyle\frac{\beta x_1x_5}{\kappa+x_5}+\omega x_4-\mu x_1 \right)\\
&\quad + \lambda_2 \left( \displaystyle\frac{\beta y_1y_5}{\kappa+y_5}
-(\delta u+\alpha_1+\mu)x_2 \right)+ \lambda_3 \Big (\delta 
ux_2-(\varepsilon+\alpha_2+\mu)x_3 \Big )\\
&\quad + \lambda_4 \Big ( \varepsilon x_3-(\omega+\mu)x_4 \Big ) 
+ \lambda_5 (\eta x_2-dx_5).
\end{split}
\end{equation}
Let the pair $\big (X^\diamond(\cdot),u^\diamond(\cdot)\big )
\in\mathcal{X}\times\mathcal{U}$ be an optimal solution of \eqref{mincostfunct}. 
Then, in agreement with \cite{Gollmann1,Gollmann2}, 
there is an adjoint function $\lambda^\diamond(\cdot)=\big(\lambda_1^\diamond(\cdot),
\lambda_2^\diamond(\cdot),\lambda_3^\diamond(\cdot),\lambda_4^\diamond(\cdot),
\lambda_5^\diamond(\cdot)\big ) \in W^{1,1}\big ([0,T];\mathbb{R}^5\big )$ such that 
the following conditions hold a.e. in $t \in [0,T]$:
\begin{enumerate}
\item[1)] the differential adjoint equations:
\begin{equation}
\label{CN_dif_adj_system}
\begin{cases}
\begin{split}
\dot{\lambda_1^\diamond}(t)
&=-\frac{\partial H}{\partial x_1}[t]
-\frac{\partial H}{\partial y_1}[t+\tau]\chi_{[0,T-\tau]}(t)\\
&=\frac{\beta x_5^\diamond(t)}{k+x_5^\diamond(t)}\Big (\lambda_1^\diamond(t)
-\lambda_2^\diamond(t+\tau)\chi_{[0,T-\tau]}(t)\Big )+\mu\lambda_1^\diamond(t),\\[0.2 cm]
\dot{\lambda_2^\diamond}(t)
&=-\frac{\partial H}{\partial x_2}[t] 
= -W_I +{\lambda_2^\diamond(t)}\Big (\delta u^\diamond(t)
+ \alpha_1 + \mu \Big ) - \delta\lambda_3^\diamond(t) u^\diamond(t) 
- \eta\lambda_5^\diamond(t),\\[0.2 cm]
\dot{\lambda_3^\diamond}(t)
&=-\frac{\partial H}{\partial x_3}[t] 
= \Big (\varepsilon+\alpha_2+\mu\Big )\lambda_3^\diamond(t)
-\varepsilon\lambda_4^\diamond(t),\\[0.2cm]
\dot{\lambda_4^\diamond}(t)
&=-\frac{\partial H}{\partial x_4}[t] = -\omega\lambda_1^\diamond(t) 
+ \Big (\omega+\mu\Big )\lambda_4^\diamond(t),\\[0.2cm]
\dot{\lambda_5^\diamond}(t)
&=-\frac{\partial H}{\partial x_5}[t]
-\frac{\partial H}{\partial y_5}[t+\tau]\chi_{[0,T-\tau]}(t)\\
&= -W_B + \frac{\beta k x_1^\diamond(t)}{\big (k+x_5^\diamond(t)\big )^2}
\Big ( \lambda_1^\diamond(t)-\lambda_2^\diamond(t+\tau)
\chi_{[0,T-\tau]}(t)\Big )+d\lambda_5^\diamond(t);\\[0.1cm]
\end{split}
\end{cases}
\end{equation}
\item[2)] the transversality conditions 
$$
\lambda^\diamond(T)=0,
$$
in view of the free terminal state $X(T)$;
\item[3)] the minimality condition:
\begin{equation}
\label{OCP:min_nec_cond_u}
H\big(X^\diamond(t),y_1^\diamond(t),y_5^\diamond(t),u^\diamond(t),
\lambda^\diamond(t)\big) = \min_{u \in \mathcal{U}_0} H\big(X^\diamond(t),
y_1^\diamond(t),y_5^\diamond(t),u,\lambda^\diamond(t)\big),
\end{equation}
where
$$\mathcal{U}_0 =
\big \{
u \in \mathbb{R} \, | \,
1 \leq u \leq u_{\max} \big \}.
$$
\end{enumerate}
Note that the Hamiltonian \eqref{Hamiltonian} is linear 
in the control variable. Hence, the minimizing control
is determined by the sign of the switching function
\begin{equation}
\label{switching-function}
\phi(t)= \frac{\partial H}{\partial u}[t]
= W_u + \delta x_2^\diamond(t)\Big(\lambda_3^\diamond(t)-\lambda_2^\diamond(t)\Big)
\end{equation}
as
\begin{equation}
\label{control-law-q1}
u^\diamond(t)=
\begin{cases}
\begin{split}
&u_{\max}, \ &&\text{if} \quad \phi(t)<0;\\
&1,        \ &&\text{if} \quad \phi(t)>0;\\
&\text{singular},\ &&\text{if} \quad \phi(t)=0\ \forall t \in I_s\subset[0,T].
\end{split}
\end{cases}
\end{equation}
If the switching function has only finitely many isolated zeros 
in an interval $I_b \subset [0,T]$, then the
control $u^\diamond $ is called \textit{bang-bang} on $I_b$. 
In case of a \textit{singular control},
where $\phi(t) = 0$ on $I_s \subset [0,T]$, 
the Pontryagin Minimum Principle does not directly provide values of the control.
We shall not discuss the singular case here, 
since in our computations we never 
encountered \textit{singular controls}.

\begin{rmrk}
When we write $H[t]$, we mean that the arguments of the Hamiltonian $H$ 
are applied in $t$ along the extremal, that is,
$H[t] = H\big (X^\diamond(t),y_1^\diamond(t),y_5^\diamond(t),u^\diamond(t),
\lambda^\diamond(t)\big)$.
\end{rmrk}


\subsection{Numerical solutions of the time-delayed SIQRB control model}
\label{subsec:OC_SIQRBB_delayed}

With the purpose to obtain strategies that could have stopped 
the spread of the cholera outbreak in Haiti mentioned 
in Section~\ref{sec:num:simu}, we consider here that 
control measures are taken. The absence of a control measure 
means that we are considering $u(t)=1$ for all time $t\in[0,T]$, i.e., 
all infected individuals are moved to quarantine at the end of 50 days. 
Our aim consists in obtaining a faster response with the minimum cost. 
To illustrate this, we are going to solve optimal control problem 
\eqref{mincostfunct}, considering several cases. For solving it, we apply
the discretization methods developed in \cite{Gollmann1,Gollmann2}.
The resulting large-scale non-linear programming problem (\textsc{NLP}) 
can be conveniently formulated using the Applied Modeling Programming Language 
-- \textsc{AMPL} (see \cite{AMPL-0,AMPL}), which can be linked to several 
efficient optimization solvers. We use the 
Interior-Point optimization solver -- \textsc{IPOPT}, 
developed by W\"achter and Biegler \cite{Waechter-Biegler}, 
with tolerance $tol = 10^{-8}$.
As integration method, we implement Euler's method with $T=182$ days 
and at least $12740$ grid points. In the following, we shall compare 
several solutions by considering three different combinations 
of weights $W_I$ and $W_B$:
\begin{enumerate}
\item[1)] Case 1: $W_I = W_B = 1$;
\item[2)] Case 2: $W_I = 10$ and $W_B = 1$;
\item[3)] Case 3: $W_I = 1$ and $W_B = 10$.
\end{enumerate}
For these computations, we use the other parameters values 
and the initial conditions of Table~\ref{Tab_Parameter}.

\begin{figure}[ht!]
\centering
\includegraphics[scale=0.55]{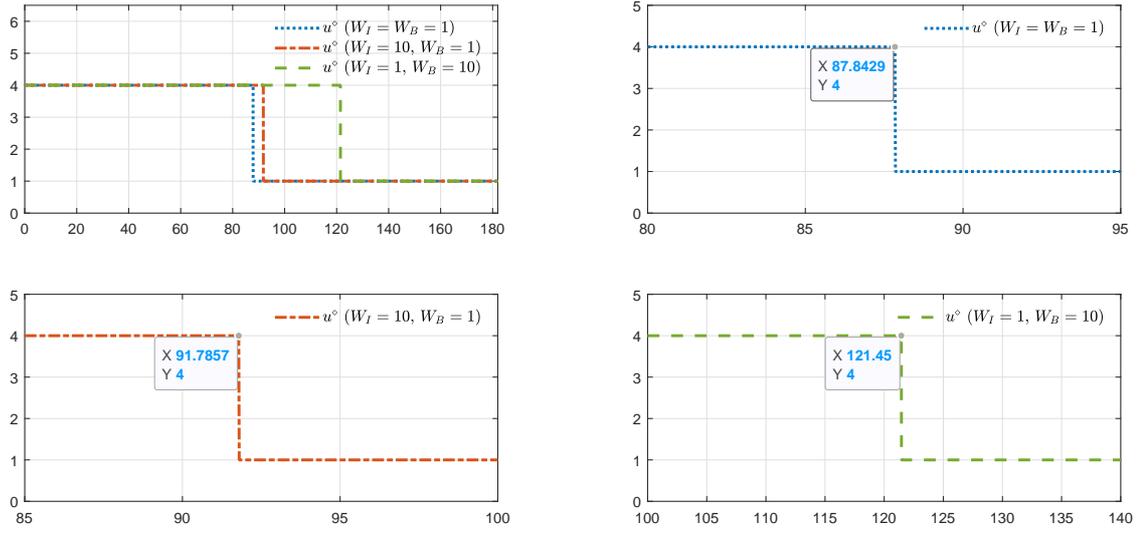}
\caption{Extremal control $u^\diamond$ and zoom of it for Cases~1, 2 and 3 
(dashed blue, dashed orange and dashed light green curves, respectively) associated 
with the optimal control problem \eqref{mincostfunct}, 
considering parameter values of Table~\ref{Tab_Parameter}.}
\label{Fig-only-u}
\end{figure}

\begin{figure}[ht!]
\centering
\includegraphics[scale=0.5]{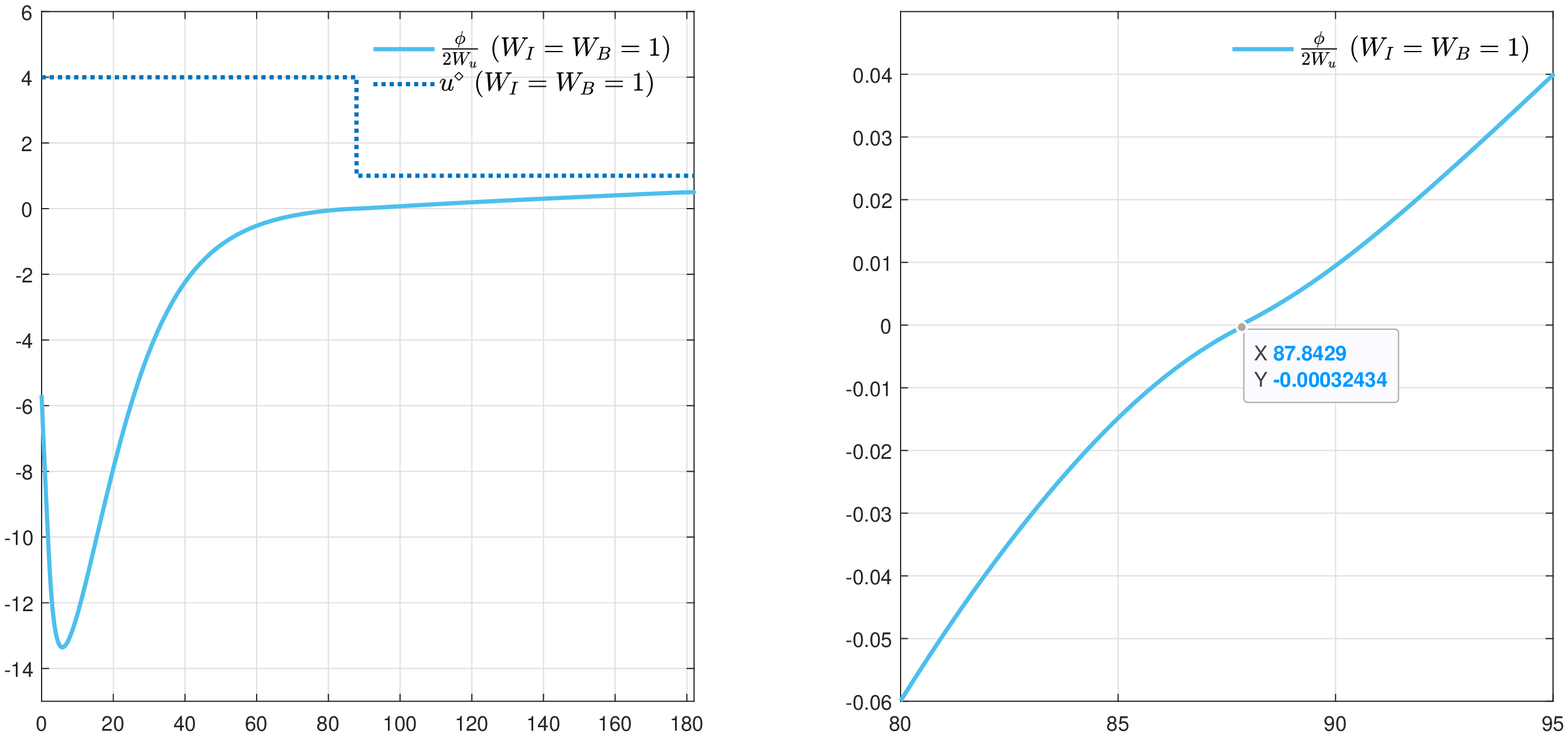}
\caption{Extremal control $u^\diamond$, 
satisfying control law \eqref{control-law-q1},	
and $\frac{\phi}{2W_u}$ (left column), 
as well as a zoom of $\frac{\phi}{2W_u}$ (right column), 
for problem \eqref{mincostfunct} in Case~1, 
with the parameter values of Table~\ref{Tab_Parameter}.}
\label{Fig-phi-case1}
\end{figure}

\bigskip
\noindent \textbf{Case 1 ($\mathbf{W_I = W_B = 1}$}). 
The obtained cost functional value and initial values 
of the adjoint variables are, approximately, given by
$$
J_1(X^\diamond,u^\diamond) \simeq 4.299001 
\times 10^6 \quad \text{and} \quad
\lambda^\diamond(0) 
\simeq (363.76,\ 379.71,\ 15.875,\ 17.704,\ 3.1978).
$$
The control is \textit{bang-bang} with one switch at $t = t^s$:
\begin{equation}
\label{control-q=1-tau=0}
u^\diamond (t) =
\begin{cases}
\begin{split}
& 4, && \mbox{if} \quad  t \in \left[0,t^s\right], \\
& 1, && \mbox{if} \quad  t \in \left]t^s,T\right],
\end{split}
\end{cases}
\end{equation}
where 
$t^s \simeq 87.843$
(see top left and right column of Figure~\ref{Fig-only-u}). 
Note that the switching function satisfies the strict \textit{bang-bang property} 
(see right column of Figure~\ref{Fig-phi-case1}):
\begin{equation}
\label{strict:bang:bang:prop:case1}
\begin{cases}
\phi(t)<0 \text{ for } t \in \left[0,t^s\right[,\\
\dot{\phi}(t^s)>0,\\
\phi(t)>0 \text{ for } t \in \left]t^s,T\right].
\end{cases}
\end{equation}

\begin{figure}[ht!]
\centering
\includegraphics[scale=0.5]{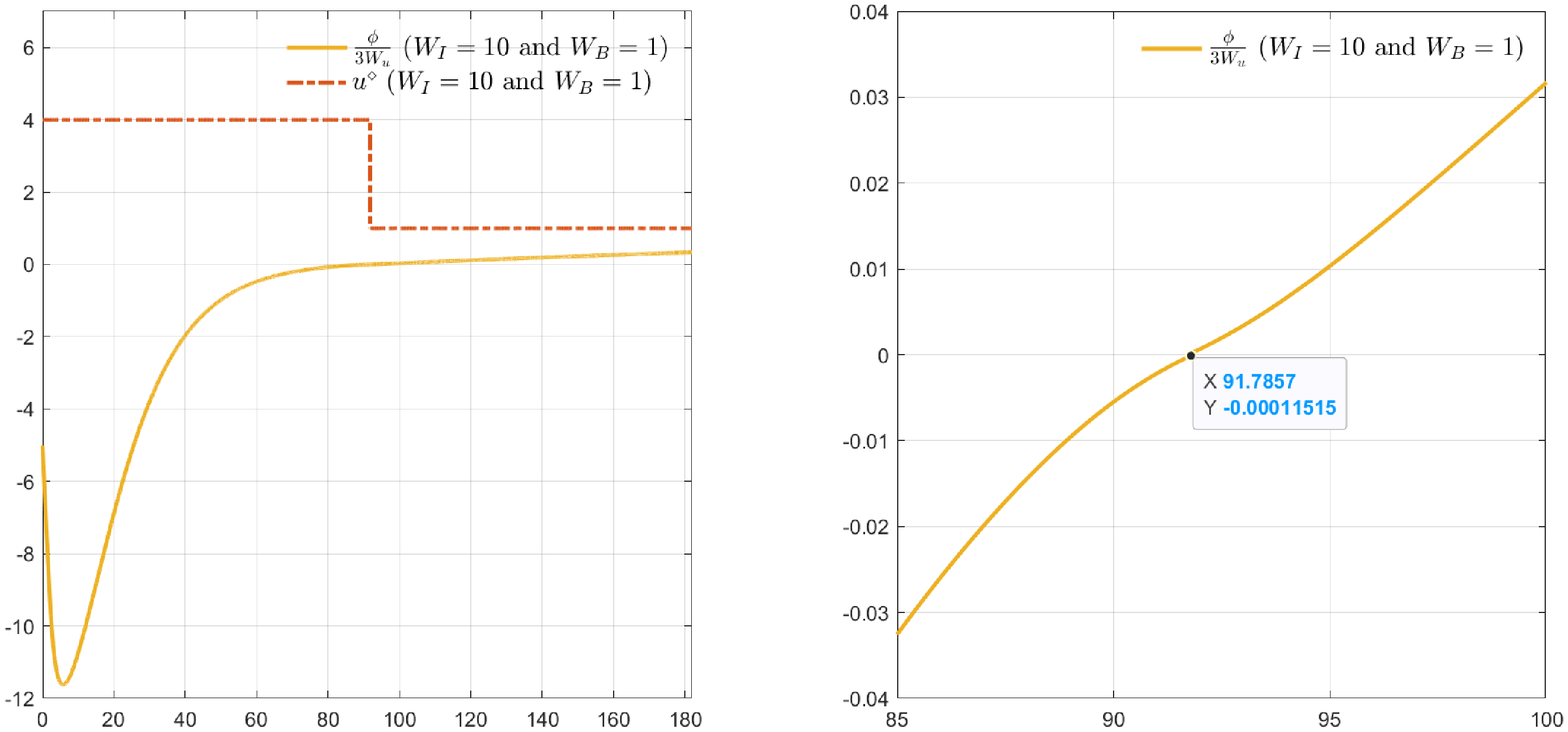}
\caption{Extremal control $u^\diamond$, 
satisfying control law \eqref{control-law-q1},	
and $\frac{\phi}{3W_u}$ (left column), 
as well as a zoom of $\frac{\phi}{3W_u}$ (right column), 
for problem \eqref{mincostfunct} in Case~2,
with the parameter values of Table~\ref{Tab_Parameter}.}
\label{Fig-phi-case2}
\end{figure}

\bigskip
\noindent \textbf{Case 2 ($\mathbf{W_I = 10}$ and $\mathbf{W_B = 1}$)}. 
The obtained cost functional value and initial values 
of the adjoint variables are, approximately, given by
$$
J(X^\diamond,u^\diamond) \simeq 5.166139 \times 10^6
\quad \text{and} \quad
\lambda^\diamond(0) \simeq (465.69,\ 490.35,\ 18.875,\ 21.210,\ 3.2784).
$$
The control has also the \textit{bang-bang} structure \eqref{control-q=1-tau=0}
with $t^s \simeq 91.79$
(see left column of Figure~\ref{Fig-only-u}) and 
the switching function also satisfies the strict \textit{bang-bang property}, 
as in \eqref{strict:bang:bang:prop:case1} 
(see right column of Figure~\ref{Fig-phi-case2}).

\begin{figure}[ht!]
\centering
\includegraphics[scale=0.5]{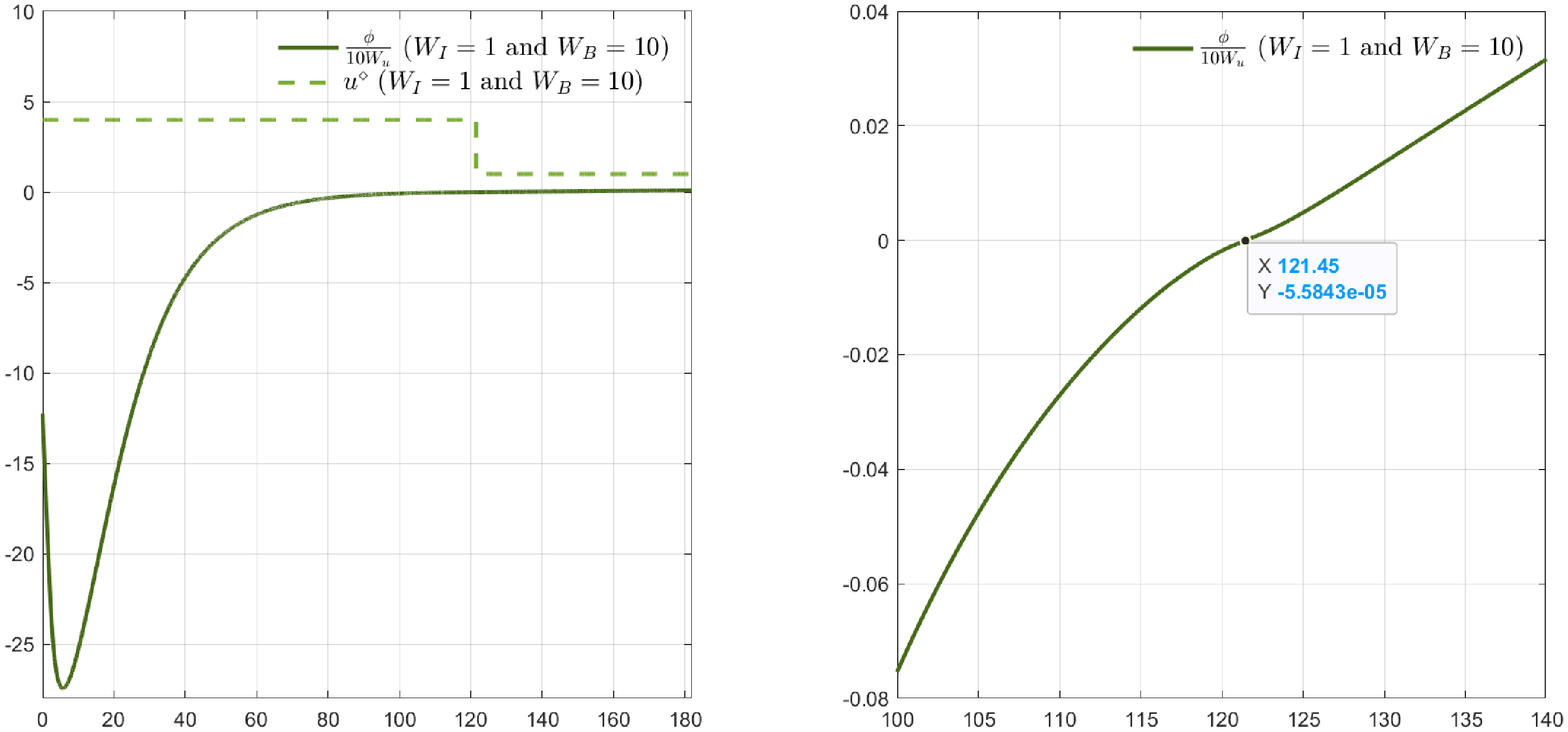}
\caption{Extremal control $u^\diamond$, 
satisfying control law \eqref{control-law-q1}, 	
and $\frac{\phi}{10W_2}$ (left column), 
as well as a zoom of $\frac{\phi}{10W_u}$ (right column), 
for problem \eqref{mincostfunct} in Case~3,
with the parameter values of Table~\ref{Tab_Parameter}.}
\label{Fig-phi-case3}
\end{figure}

\bigskip 
\noindent \textbf{Case 3 ($\mathbf{W_I=1}$ and $\mathbf{W_B=10}$)}. 
The obtained cost functional value and initial values of the 
adjoint variables are, approximately, given by
$$
J(X^\diamond,u^\diamond) \simeq 3.7821542 \times 10^7
\quad \text{and} \quad
\lambda^\diamond(0) \simeq (3450.2,\ 3734.6,\ 100.66,\ 117.78,\ 32.964).
$$
Again, the control has the \textit{bang-bang} 
structure \eqref{control-q=1-tau=0} with $t^s \simeq 121.45$
(see top left and bottom right column of Figure~\ref{Fig-only-u})
and the switching function satisfies the strict \textit{bang-bang property} 
(see right column of Figure~\ref{Fig-phi-case3}), as in \eqref{strict:bang:bang:prop:case1}.

\begin{figure}[ht!]
\centering
\includegraphics[scale=0.5]{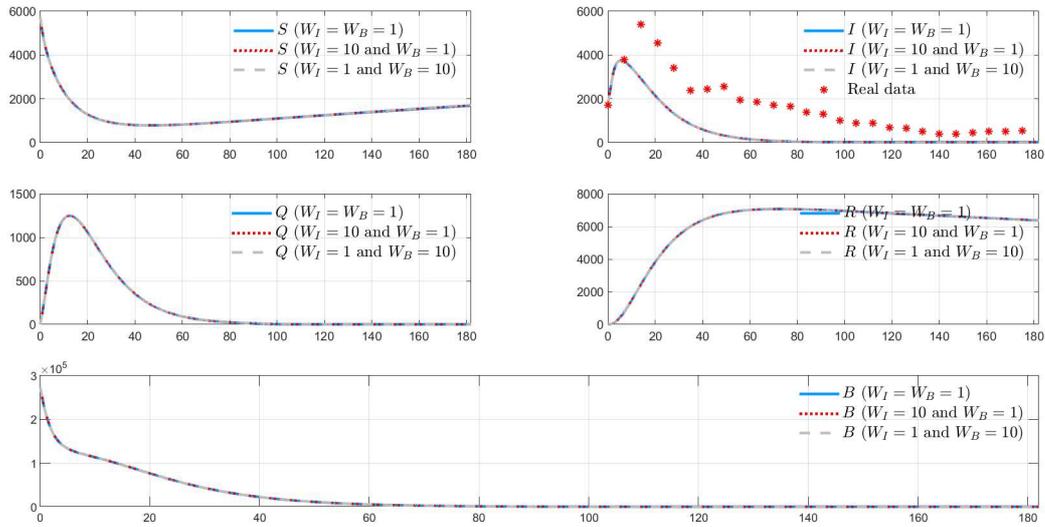}
\caption{Extremal trajectories $S^\diamond$, $I^\diamond$, $Q^\diamond$, 
$R^\diamond$ and $B^\diamond$ of the optimal control problem \eqref{mincostfunct}, 
considering parameter values of Table~\ref{Tab_Parameter}, for Cases 1, 2 and 3 
(solid light blue, dashed red, and dashed grey curves, respectively).}
\label{Fig-SIQRB-comparison}
\end{figure}

Figures~\ref{Fig-only-u}--\ref{Fig-phi-case3} summarize, in the three cases, 
our numerical findings associated with the extremal controls, 
while Figure~\ref{Fig-SIQRB-comparison} represents 
the corresponding extremal trajectories. 
The main message of
Figures~\ref{Fig-phi-case1}--\ref{Fig-phi-case3} is that the  
controls we have obtained numerically satisfy, always, 
the control law \eqref{control-law-q1}. Hence, the necessary 
conditions are satisfied and thus we have found the \emph{extremals}.
However, we remark that it does not suffice to check the 
strict bang-bang property to conclude that extremal controls are optimal. 
For time-delayed optimal control problems with bang-bang solutions, 
it remains an open question to prove a second order sufficient optimality condition.

\begin{rmrk}
With respect to all considered cases, we present in Figures~\ref{Fig-phi-case1}, 
\ref{Fig-phi-case2} and \ref{Fig-phi-case3} function $\frac{\phi}{cW}$, where 
$c$ is a positive constant, instead of $\phi$, because the values obtained by $\phi$ 
are much bigger than those achieved by control $u^\diamond$. Thus, using 
function $\frac{\phi}{cW}$, it is possible to compare the signal of $\phi$ 
with the behaviour of the extremal control $u^\diamond$ in the same plot.
\end{rmrk}

\begin{figure}[ht!]
\centering
\includegraphics[scale=0.5]{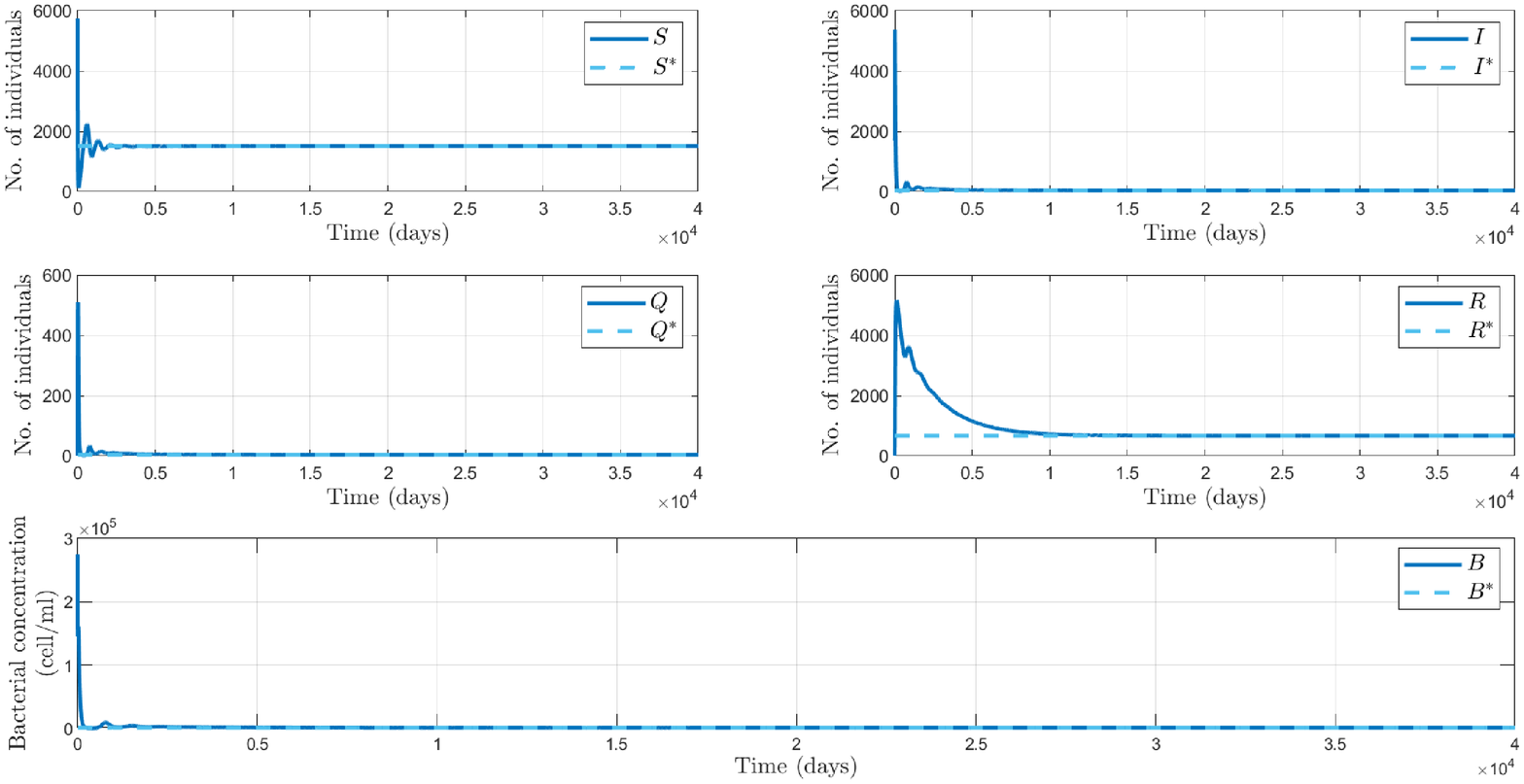} 
\caption{Asymptotic stability of the endemic equilibrium $(S^*,I^*,Q^*,R^*,B^*)$ 
of model \eqref{ModeloColera_delay}, using all the values of Table~\ref{Tab_Parameter}.}
\label{FigHaitiComQ_delayed_all_solutions_stability}
\end{figure}

Figure~\ref{Fig-SIQRB-comparison} shows the comparison 
of the extremal state trajectories $S^\diamond$, $I^\diamond$, 
$Q^\diamond$, $R^\diamond$ and $B^\diamond$ for the
optimal control problem \eqref{mincostfunct}, 
considering the parameter values of Table~\ref{Tab_Parameter}, 
for the three studied cases. 
Note that, although one can see that the obtained trajectories
are very similar in all the cases,
the controls corresponding to the three cases are very different,
as well as the associated costs (see Figure~\ref{Fig-only-u})
which implies different policies to control the outbreak.
Furthermore, we also observed numerically that the endemic equilibrium 
\eqref{EndemicEquilibrium} of the delayed model \eqref{ModeloColera_delay} 
is locally asymptotic stable, when we consider all the values 
of Table~\ref{Tab_Parameter} 
(see Figure~\ref{FigHaitiComQ_delayed_all_solutions_stability}).


\section{Conclusions}
\label{sec:conclusion}

The paper is devoted to the formulation and analysis of a time-delayed mathematical model
for cholera. We have improved our recent results by adding a time delay, which represents 
the time between the instant at which an individual becomes infected and the instant
at which he begins to have symptoms of cholera infection by the bacterium \emph{Vibrio cholerae}.
The considered model is a SIQR (Susceptible--Infectious--Quarantined--Recovered) system, 
where an additional class $B$ is considered: a class of bacterial concentration
in the dynamics of cholera. The host population is divided into four classes: susceptible, 
infectious with symptoms, quarantined, and recovered. An additional compartment
reflects the bacterial concentration. The formulated model is analyzed, providing the
non-negativity of the solutions for non-negative initial conditions, as well as the disease-free 
equilibrium, basic reproduction number, and endemic equilibrium. For positive time
delays, the stability of the equilibrium points is also analyzed. Then, we consider
the cholera outbreak that occurred in Haiti in 2010 and 2011. Finally, a SIQRB time-delayed 
optimal control problem, with linear cost functional, is formulated
and analyzed numerically.

The delayed model is more realistic and describes 
better the reality, since the symptoms of cholera disease can appear 
from a few hours up to five days after infection.
Usually, the symptoms appear in 2--3 days after infection (see \cite{TimeDelayCholera}).
We proved the positivity of the solutions of the delayed SIQRB model and analyzed 
the local asymptotic stability of the equilibrium points. Our analysis shows 
that the ingestion rate of the bacteria through contaminated sources $\beta$ 
has an important influence on the stability of the endemic equilibrium.
We improved the numerical simulations done in \cite{Lemos-Paiao} for the cholera 
outbreak in the Department of Artibonite -- Haiti, from 1 November 2010 to 
1 May 2011 \cite{WHO}, and showed that the delayed model fits better 
this cholera outbreak. A delayed optimal control problem was proposed and analyzed
with a $L^1$ cost functional, where the objective is to minimize 
the number of infective individuals and the bacterial concentration, as well as 
the cost of interventions associated with the control treatment through quarantine. 
Necessary optimality conditions were applied to all considered cases of the optimal 
control problem. The delayed optimal control problems were solved numerically, 
applying the discretization methods developed in \cite{Gollmann1,Gollmann2}.
Summarizing, the main contributions and novelties of our paper are: 
(i) the proposal of a more realistic cholera model with time delays;
(ii) a local asymptotic stability analysis with respect to the proposed cholera model;
(iii) the proposal, for the new model, of a more appropriate optimal control 
problem, in the biomedical framework, with linear cost measure,
instead of a quadratic one as considered before in the literature;
(iv) a numerical study of the proposed model and optimal control problem in order to 
translate cholera Haiti outbreak and to provide a control measure 
to stop the disease spread, respectively.

As future research directions, it remains to obtain 
(i) sufficient conditions for the local asymptotic stability of the endemic equilibrium, 
as well as the study of global stability, and (ii) second order sufficient optimality 
conditions for time-delayed optimal control problems with bang-bang solutions.
Furthermore, as we have some parameters for which we only know an interval of values 
where they belong, it could be interesting to relate some of them under inequality constraints. 
Thus, we also point out, as future work, the study of parametric optimal control problems
for cholera, by applying parameter sensitivity analysis.


\begin{acknowledgement}
This research was supported 
by the Portuguese Foundation for Science and Technology (FCT)
within projects UIDB/04106/2020 and UIDP/04106/2020 (CIDMA). 
Lemos-Pai\~{a}o was also supported by the Ph.D. 
fellowship PD/BD/114184/2016 and  
Silva by the FCT Researcher Program CEEC Individual 2018
with reference CEECIND/00564/2018. The authors are grateful 
to an Associate Editor and two anonymous Referees, 
who kindly reviewed an earlier version of the manuscript 
and provided several valuable suggestions and comments.
\end{acknowledgement}


\clearpage



\end{document}